\documentclass[10pt,a4paper]{article}
\usepackage{cite}
\usepackage{amscd}
\usepackage[latin1]{inputenc}
\usepackage{amsmath}
\usepackage{amsfonts}
\usepackage{amssymb}
\usepackage{color}
\usepackage{float}
\usepackage{amsthm}
\usepackage{graphicx}
\usepackage{listings}
\newtheorem{theorem}{Theorem}

\newtheorem{proposition}[theorem]{Proposition}

\def\Qed{\hfill\raisebox{.6ex}{\framebox[2.5mm]{}}\\[.15in]}

\def\m{\mathbb}
\def\mc{\mathcal}

\begin{document}
\date{}
\title{New surfaces with $K^2=7$ and $p_g=q\leq 2$}
\author{Carlos Rito}
\maketitle

\begin{abstract}

We construct smooth minimal complex surfaces of general type with $K^2=7$ and:
$p_g=q=2,$ Albanese map of degree $2$ onto a $(1,2)$-polarized abelian surface;
$p_g=q=1$ as a double cover of a quartic Kummer surface;
$p_g=q=0$ as a double cover of a numerical Campedelli surface with $5$ nodes.

\noindent 2010 MSC: 14J29.

\end{abstract}

\section{Introduction}

Despite the efforts of several authors in past years, surfaces of general type with the lowest possible value
of the holomorphic Euler characteristic $\chi=1$ are still not classified. For these surfaces the geometric
genus $p_g,$ the irregularity $q$ and the self-intersection of the canonical divisor $K$ satisfy:
$$1+p_g\leq K^2\leq 9\ \ \ {\rm if}\ \ \ p_g\leq 1,$$
$$\ \ \ \ \ \ \ 4\leq K^2\leq 9\ \ \ {\rm if}\ \ \ p_g=2,$$
$$\ \ \ \ \ \ K^2=6\ {\rm or}\ 8\ \ \ {\rm if}\ \ \ p_g=3,$$
$$\ \ \ \ \ \ \ \ \ \ \ \ K^2=8\ \ \ {\rm if}\ \ \ p_g=4,$$
from the Bogomolov-Miyaoka-Yau and Debarre inequalities, \cite{HP} and the Beauville Appendix in \cite{De}
(cf. also \cite{CCM}, \cite{Pi}).

According to Sai-Kee Yeung, the case with $p_g=q=2$ and $K^2=9$ does not occur (see Section 6 in the revised
version of the paper \cite{Ye}, available at\newline 
\verb%http://www.math.purdue.edu/~yeung/%).

So there are examples for all possible values of the invariants except for one mysterious case:
$$K^2=7, p_g=q=2.$$ The cases $K^2=7,$ $p_g=q=1$ or $0$ are also intriguing:
\begin{itemize}
\item $p_g=q=1.$ Lei Zhang \cite{Zh} has shown that one of three cases occur:\\ a) the bicanonical map is
birational;\\ b) the bicanonical map is of degree $2$ onto a rational surface;\\ c) the bicanonical map is
of degree $2$ onto a Kummer surface.\\ The author has given examples for a) \cite{Ri2} and b) \cite{Ri3},
but so far it is not known if c) can occur.
\item $p_g=q=0.$ Yifan Chen \cite{Ch2} considers the case when the automorphism group of the surface $S$
contains a subgroup isomorphic to $\mathbb Z_2^2.$ He shows that three different families of surfaces may exist:\\
a) $S$ is an Inoue surface \cite{In};\\
b) $S$ belongs to the family constructed by him in \cite{Ch1};\\
c) a third case, in particular $S$ is a double cover of a surface with $p_g=q=0$ and
$K^2=2$ with $5$ nodes.\\
The existence of this last case is an open problem.
\end{itemize}

In this paper we show the existence of the above three open cases.
We give constructions for surfaces with $K^2=7$ and:
\begin{itemize}
\item $p_g=q=2,$ Albanese map of degree $2$ onto a $(1,2)$-polarized abelian surface;
\item $p_g=q=1,$ bicanonical map of degree $2$ onto a Kummer surface;
\item $p_g=q=0$ as a double cover of a numerical Campedelli surface with $5$ nodes.
\end{itemize}
In all cases the surface can be seen as a double cover with branch locus as in the result below.
In particular we show that a construction for the case $p_g=q=2$ as suggested by Penegini and Polizzi
\cite[Remark 2.2]{PP} does exist.

\begin{proposition}\label{prop1}
Let $X$ be an Abelian, $K3$ or Enriques surface containing $n$ disjoint
$(-2)$-curves $A_1,\ldots,A_n$, $n=0, 16$ or $8,$ respectively.
Assume that $X$ contains a reduced curve $B$ and a divisor $L$ such that $$B+\sum_1^n A_i\equiv 2L,$$
$B$ is disjoint from $\sum_1^n A_i,$
$B^2=16$ and $B$ contains a $(3,3)$-point and no other singularity.
Let $S$ be the smooth minimal model of the double cover of $X$ with branch locus $B+\sum_1^n A_i.$
Then $\chi(\mathcal O_S)=1$ and $K_S^2=7.$
\end{proposition}
\noindent{\bf Proof :}\\
This follows from the double cover formulas (see e.g. \cite[V.22]{BHPV}) and the fact that
a $(3,3)$-point decreases both $\chi$ and $K^2$ by $1$ (see e.g. \cite[p. 185]{Pe}):
$$\chi(\mathcal O_S)=2\chi(\mathcal O_X)+\frac{1}{2}L(K_X+L)-1=1,$$
$$K_S^2=2(K_X+L)^2+n-1=7.$$\Qed

\bigskip
\noindent{\bf Notation}

We work over the complex numbers. All varieties are assumed to be projective algebraic.
A $(-n)$-curve on a surface is a curve isomorphic to $\m P^1$ with self-intersection $-n.$
An $(m_1,m_2)$-point of a curve, or point of type $(m_1,m_2),$ is a singular point of multiplicity
$m_1$ which resolves to a point of multiplicity $m_2$ after one blow-up.
Linear equivalence of divisors is denoted by $\equiv.$
The rest of the notation is standard in Algebraic Geometry.\\

\bigskip
\noindent{\bf Acknowledgements}

The author is a member of the Center for Mathematics of the University of Porto.
This research was partially supported by FCT (Portugal) under the project PTDC/MAT-GEO/0675/2012 and by CMUP (UID/MAT/00144/2013), which is funded by FCT with national (MEC) and European structural funds through the programs FEDER, under the partnership agreement PT2020.

\section{Bidouble covers}
A bidouble cover is a finite flat Galois morphism with Galois group $\mathbb Z_2^2.$
Following \cite{Ca} or \cite{Pa}, to define a bidouble cover $\pi:V\rightarrow X,$
with $V,$ $X$ smooth surfaces, it suffices to present:
\begin{itemize}
\item smooth effective divisors $D_1, D_2, D_3\subset X$ with pairwise
transverse intersections and no common intersection;
\item line bundles $L_1, L_2, L_3$ such that $L_g+D_g\equiv L_j+L_k$ for
each permutation $(g,j,k)$ of $(1,2,3).$
\end{itemize}
%
%
One has
$$\chi(\mc O_V)=4\chi(\mc O_X)+\frac{1}{2}\sum_1^3 L_i(K_X+L_i),$$
$$p_g(V)=p_g(X)+\sum_1^3 h^0(X,\mc O_X(K_X+L_i))$$
and $$2K_V\equiv\pi^*\left(2K_X+\sum_1^3 L_i\right),$$
which implies $$K_V^2=\left(2K_X+\sum_1^3 L_i\right)^2.$$

\section{Example with $p_g=q=2$}\label{section3}

\noindent{\bf Step 1}\\
Let $T_1,\ldots,T_4\subset\m P^2$ be distinct lines through a point $p_0,$ let $p_1,p_2\ne p_0$ be
points in $T_1,T_2,$ respectively, and $C_1,C_2$ be distinct smooth conics tangent to $T_1,T_2$ at
$p_1,p_2.$
Consider the map $$\mu:X\longrightarrow\m P^2$$ which resolves the singularities of the divisor
$C_1+C_2+T_1+\cdots+T_4.$ Then $\mu$ is given by blow-ups at
$$p_0,p_1,p_1',p_2,p_2',p_3,\ldots,p_{10},$$ where $p_i'$ is the point infinitely near to $p_i$
corresponding to the line $T_i,$
and $p_3,\ldots,p_{10}$ are nodes of $C_1+C_2+T_3+T_4.$
Let $E_0,E_1,E_1',E_2,E_2',E_3,\ldots,E_{10}$ be the corresponding exceptional divisors
(with self-intersection $-1$) and let $$\pi:V\longrightarrow X$$ be the bidouble cover defined by the divisors
$$
\begin{array}{l}
D_1:=\left(\widetilde{T_1}+\widetilde{T_2}-2E_0-2E_1'-2E_2'\right)+\sum_3^{10}E_i,\\
D_2:=\widetilde{T_3}+\widetilde{T_4}-2E_0-\sum_3^{10}E_i,\\
D_3:=\widetilde{C_1}+\widetilde{C_2}-2E_1-2E_1'-2E_2-2E_2'-\sum_3^{10}E_i,
\end{array}
$$
where the notation $\widetilde{\cdot}$ stands for the total transform $\mu^*(\cdot).$

Notice that $D_1$ is the union of $\sum_3^{10}E_i$ with four $(-2)$-curves contained in the pullback of $T_1+T_2,$
and $D_2,D_3$ are just the strict transforms of $T_3+T_4,$ $C_1+C_2,$ respectively.

One can easily see that the divisors $D_1,$ $D_2$ and $D_3$ have pairwise transverse intersections and no common intersection.

Denote by $T$ a general line of $\m P^2$ and let
$$
\begin{array}{l}
L_1:=3\widetilde{T}-E_0-E_1-E_1'-E_2-E_2'-\sum_3^{10}E_i,\\
L_2:=3\widetilde{T}-E_0-E_1-2E_1'-E_2-2E_2',\\
L_3:=2\widetilde{T}-2E_0-E_1'-E_2'.
\end{array}
$$
Then
$$
\begin{array}{l}
K_X+L_1\equiv 0,\\
K_X+L_2\equiv-E_1'-E_2'+\sum_3^{10}E_i,\\
K_X+L_3\equiv-\widetilde{T}-E_0+E_1+E_2+\sum_3^{10}E_i
\end{array}
$$
and
$$\chi(\mc O_V)=4+\frac{1}{2}(0-4-4)=0,$$ $$p_g(V)=0+1+0+0=1.$$

Let $X_1$ be the surface given by the double covering $\phi:X_1\longrightarrow X$ with branch locus $D_2+D_3.$
The divisor $\phi^*\left(\widetilde{T_1}+\widetilde{T_2}-2E_0-2E_1'-2E_2'\right)$ is a disjoint union of $8$ $(-2)$-curves,
and the divisor $\phi^*\left(\sum_3^{10}E_i\right)$ is also a disjoint union of $8$ $(-2)$-curves.
Hence $\phi^*(D_1)$ is a disjoint union of $16$ $(-2)$-curves.
The canonical divisor $K_V$ of $V$ is the support of the pullback of $D_1,$ a disjoint union of
$16$ $(-1)$-curves. So the minimal model $V'$ of $V$ is an abelian surface, with Kummer surface $X_1.$
Notice that the lines $T_1,\ldots,T_4$ give rise to elliptic fibres of type $I_0^*$ in $X_1$
(four disjoint $(-2)$-curves plus an elliptic curve with multiplicity $2$).\\

\noindent{\bf Step 2}\\
Now let $R$ be the tangent line to $C_1$ at $p_3\in C_1\cap T_3.$
We claim that the strict transform $\widehat{R}\subset V'$ of $R$ is a curve with a tacnode
(singularity of type $(2,2)$) at the pullback of $p_3$ and with self-intersection $\widehat R^2=8.$
In fact, the covering $\pi$ factors as
$$V\xrightarrow{\ \ \ \varphi\ \ \ } X_1 \xrightarrow{\ \ \ \phi\ \ \ } X.$$
The strict transforms $R', C_1'\subset X$ of $R, C_1$ meet at a point in the $(-1)$-curve $E_3.$
Since $C_1'$ is contained in the branch locus of the covering $\phi,$ then the curve $R'':=\phi^*(R')$
is tangent to the $(-2)$-curve $\phi^*(E_3).$ This curve is in the branch locus of $\varphi,$ hence
the curve $R''':=\varphi^*(R'')$ has a node at a point $p$ in the $(-1)$-curve
$$\overline E_3:=\frac{1}{2}(\phi\circ\varphi)^*(E_3).$$
So the image  of $R'''$ in the minimal model $V'$ of $V$ is a curve $\widehat R$ with a tacnode.

\begin{figure}[ht]
\centering
\includegraphics[width=12cm]{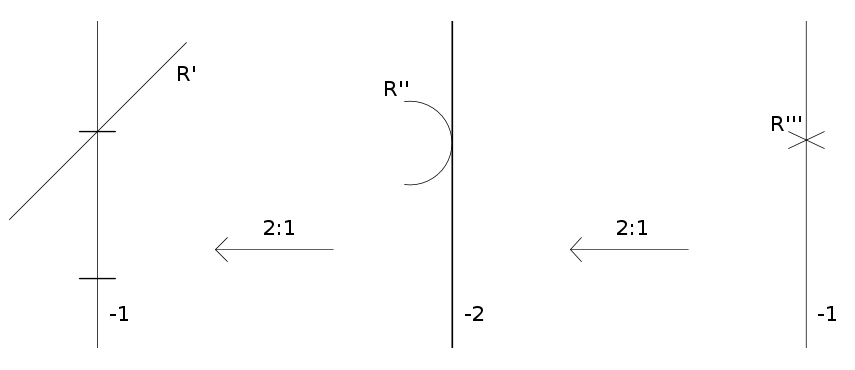}
\vspace*{-0.2cm}
\end{figure}

The reduced strict transform of the conic $C_1$ passes through $p,$ hence its image
$\widehat C_1\subset V'$ is tangent to $\widehat R$ at the tacnode.
So the divisor $\widehat R+\widehat C_1$ is reduced and has a singularity of type $(3,3).$
We want to show that it is even, i.e. there is a divisor $L$ such that
$$\widehat R+\widehat C_1\equiv 2L,$$
and that $$(\widehat R+\widehat C_1)^2=16.$$\\

\noindent{\bf Step 3}\\
The pencil of lines through the point $p_0$ induces an elliptic fibration of the surface $V.$
For $i=1,2,$ the line $T_i$ gives rise to a fibre (counted twice) which is the union of disjoint $(-1)$-curves
$\xi_1^i,\ldots,\xi_4^i$ with an elliptic curve $T_i'$ such that $\xi_j^i T_i'=1.$
These curves can be labeled such that $\xi_1^i,\xi_2^i$ correspond to the strict transform of $T_i$
and $\xi_3^i,\xi_4^i$ correspond to the $(-2)$-curve $E_i-E_i'.$
The curve $R'''$ meets $\xi_1^1,\xi_2^1,\xi_1^2,\xi_2^2,$ thus $\widehat R^2=4+4=8$
and then $(\widehat R+\widehat C_1)^2=8+0+2\times 4=16.$

Let $H$ be the line through the points $p_1$ and $p_2.$
We have
$$\pi^*(R+H)=R'''+H'+2\overline E_3+\sum_1^2 (T_i'+2\xi_3^i+2\xi_4^i),$$
where $H'\subset V$ is the strict transform of $H.$ Denote by $\widehat R,$ $\widehat H,$ $\widehat T_1$
and $\widehat T_2$ the projections of $R''',$ $H',$ $T_1'$ and $T_2'$ into the minimal model $V'$ of $V.$
Then there is a divisor $L'$ such that
$$\widehat R+\widehat H+\widehat T_1+\widehat T_2\equiv 2L'.$$

The pencil of conics tangent to the lines $T_1,T_2$ at $p_1,p_2$ induces another elliptic fibration
of the surface $V'.$ The curves $\widehat C_1$ and $\widehat H$ are fibres of this fibration.
We have 
$$\widehat R+\widehat C_1+\widehat H+\widehat C_1+\widehat T_1+\widehat T_2\equiv 2(L'+\widehat C_1).$$
Since the above fibrations have elliptic bases, the sums $\widehat H+\widehat C_1$ and
$\widehat T_1+\widehat T_2$ are even, thus there exists a divisor $L$ such that
$\widehat R+\widehat C_1\equiv 2L.$\\

\noindent{\bf Step 4}\\
Finally, consider the double cover $$\rho:S'\longrightarrow V'$$ with branch locus $\widehat R+\widehat C_1,$
determined by $L.$
It follows from Proposition \ref{prop1} that the smooth minimal model $S$ of $S'$ is a surface of general type with
$\chi=1$ and $K^2=7.$ It is known that there is no smooth minimal surface of general type with
$\chi=1,$ $K^2=7$ and $q>2$ (see \cite{HP} and the Beauville Appendix in \cite{De}).
Since $q(S)\geq q(V')=2,$ we conclude that $p_g(S)=q(S)=2.$\\

Recall that $p_3\in C_1\cap T_3$ and assume that $p_4\in C_2\cap T_4.$ The branch curve $C_1+C_2+T_1+\cdots+T_4$
is determined by the points $p_0,\ldots,p_4.$ Since any two sequences of 4 points in $\mathbb P^2,$ in general position,
are projectively equivalent, we can fix $p_0,\ldots,p_3.$ This implies that our family of examples is parametrized
by a 2-dimensional open subset of $\mathbb P^2$.

\section{Example with $p_g=q=1$}\label{section4}

Let $T_1,T_2,T_3\subset\m P^2$ be distinct lines through a point $p_0$ and $p_1,p_2,p_3\ne p_0$ be
non-collinear points in $T_1,T_2,T_3,$ respectively.
For the construction of an example with $p_g=q=0$ and $K^2=7,$
Y. Chen has shown that for a general point $p_4\ne p_0,\ldots,p_3,$ there exist:
\begin{itemize}
\item an irreducible sextic curve $C_6$ with a node at $p_0,$ a tacnode at $p_i$ with tangent line $T_i,$ $i=1,2,3,$
and having a triple point at $p_4;$
\item an irreducible quintic curve $C_5$ through $p_0, p_4$ and with a tacnode at $p_i$ with tangent line $T_i,$ $i=1,2,3.$
\end{itemize}
The curves $C_5,C_6$ correspond to the curves $\tilde B_2,\tilde B_3$
given in \cite[Proposition 2.5]{Ch1}.

Let $T$ be a general line through $p_0.$
Keeping a notation analogous to the one in Section \ref{section3}, consider the map $$\mu:X\longrightarrow\m P^2$$
which resolves the singularities of the curve $C_6$ 
and let $$\pi:V\longrightarrow X$$ be the bidouble cover defined by the divisors
$$
\begin{array}{l}
D_1:=\widetilde{T}-E_0+E_4,\\
D_2:=\left(\widetilde{T_1}+\widetilde{T_2}+\widetilde{T_3}-3E_0-\sum_1^3 2E_i'\right)+\left(\widetilde{C_6}-2E_0-\sum_1^3(2E_i+2E_i')-3E_4\right),\\
D_3:=\widetilde{C_5}-E_0-\sum_1^3(2E_i+2E_i')-E_4.
\end{array}
$$
Notice that $D_2$ is the union of the strict transform of $C_6$ with six $(-2)$-curves contained in the pullback of $T_1+T_2+T_3,$
and $D_3$ is the strict transform of $C_5.$

We verify that the divisors $D_1,$ $D_2$ and $D_3$ have pairwise transverse intersections and no common intersection.
Let $\widehat C_5,\widehat C_6$ be the strict transforms of $C_5,C_6.$ These curves are disjoint from the $(-2)$-curves contained
in the pullback of $T_1+T_2+T_3.$ It is shown in \cite[Proposition 2.5]{Ch1} that the divisor $\widehat C_5+\widehat C_6+E_4$ has at most nodal
singularities. Since the line $T$ through $p_0$ is generic, the result follows.

We have
$$
\begin{array}{l}
L_1:=7\widetilde{T}-3E_0-\sum_1^3(2E_i+3E_i')-2E_4,\\
L_2:=3\widetilde{T}-E_0-\sum_1^3(E_i+E_i'),\\
L_3:=5\widetilde{T}-3E_0-\sum_1^3(E_i+2E_i')-E_4,
\end{array}
$$

$$
\begin{array}{l}
K_X+L_1\equiv 4\widetilde{T}-2E_0-\sum_1^3(E_i+2E_i')-E_4,\\
K_X+L_2\equiv E_4,\\
K_X+L_3\equiv 2\widetilde{T}-2E_0-\sum_1^3 E_i'
\end{array}
$$
and
$$2K_X+\sum_1^3 L_i\equiv 9\widetilde{T}-5E_0-\sum_1^3(2E_i+4E_i')-E_4.$$
Thus
$$\chi(\mc O_V)=4+\frac{1}{2}(-4+0-2)=1,$$
$$p_g(V)=0+0+1+0=1$$
and
$$K_V^2=-5.$$
Since the minimal model $V'$ of $V$ is obtained contracting the $12$ $(-1)$-curves contained in
$\pi^*\left(\widetilde{T_1}+\widetilde{T_2}+\widetilde{T_3}\right),$ then $K_{V'}^2=7.$

Notice that the minimal smooth resolution of the double plane $Q\rightarrow X$ with branch locus $D_1+D_3$
is a $K_3$ surface with $16$ disjoint $(-2)$-curves, and one can obtain the surface $V$ as a double cover
of $Q$ with a branch curve $B$ as in Proposition \ref{prop1}. It can be shown that the bicanonical map of $V$ factors
through this double covering. In fact, it follows from \cite{Zh} that the bicanonical map is of degree 2 onto
a Kummer surface.

Finally we can see, as in \cite[Section 3]{Ch1}, that this family of examples is parametrized
by a 3-dimensional open subset of $\mathbb P^2\times\mathbb P^1:$ the point $p_4$ moves in an open subset of
$\mathbb P^2$ and $\widetilde T-E_0$ moves in a pencil.

\section{Example with $p_g=q=0$}\label{pgq0}

In \cite[\S 4.6]{Ri2}, the author has computed points $p_0,\ldots,p_5\in\mathbb P^2$ such that there exist:
\begin{itemize}
\item an irreducible curve $C_7$ of degree $7$ with triple points at $p_0,p_5$ and tacnodes at $p_1,\ldots,p_4$
with tangent line the line $T_i$ through $p_0,p_i,$ $i=1,\ldots,4;$
\item an irreducible curve $C_6$ of degree $6$ with a node at $p_0,$ tacnodes at $p_1,\ldots,p_4$
with tangent line the line $T_i$ through $p_0,p_i,$ $i=1,\ldots,4,$
and passing through $p_5$ such that the singularity of $C_6+C_7$ at $p_5$ is ordinary.
\end{itemize}

For the readers convenience, we give in the Appendix the equations of the curves $C_6,C_7$ computed in \cite[\S 4.6]{Ri2}
(but with a different choice of $p_0,\ldots,p_5$ in order to get shorter equations) and we verify that the curves
are exactly as stated above.

We note that for generic points $p_0,\ldots,p_5\in\mathbb P^2$ there is no such curve $C_7.$
This is because the dimension of the linear system of plane curves of degree 7 is 35, and the imposition of singularities as
above puts 36 conditions. We don't know how to construct $C_7$ without using computer algebra.
Thus here we compute just one surface, and we make no considerations about the dimension of the family of surfaces.\\

Keeping a notation as above, consider the map
$$\mu:X\longrightarrow\m P^2$$ which resolves the singularities of the curve $C_7$ 
and let $$\pi:V\longrightarrow X$$ be the bidouble cover defined by the divisors
$$
\begin{array}{l}
D_1:=\left(\widetilde{T_1}-E_0-2E_1'\right)+E_5,\\
D_2:=\left(\widetilde{T_4}-E_0-2E_4'\right)+\left(\widetilde{C_6}-2E_0-\sum_1^4(2E_i+2E_i')-E_5\right),\\
D_3:=\left(\widetilde{T_2}+\widetilde{T_3}-2E_0-2E_2'-2E_3'\right)+\left(\widetilde{C_7}-3E_0-\sum_1^4(2E_i+2E_i')-3E_5\right).
\end{array}
$$
Notice that $D_1$ is the union of $E_5$ with two $(-2)$-curves contained in the pullback of $T_1,$
the divisor $D_2$ is the union of the strict transform of $C_6$ with two $(-2)$-curves contained in the pullback of $T_4,$
and $D_3$ is the union of the strict transform of $C_7$ with four $(-2)$-curves contained in the pullback of $T_2+T_3.$

To show that the divisors $D_1,$ $D_2$ and $D_3$ have pairwise transverse intersections and no common intersection,
notice that the strict transforms $\widehat C_6,\widehat C_7$ of $C_6,C_7$ meet at an unique point,
because the intersection number of $C_6$ and $C_7$ at the points $p_0,\ldots,p_5$ is $6+4\times 8 +3=41.$
It suffices to show that this point is not in $E_5.$ In the Appendix we compute that in fact the singularities of
$C_6+C_7$ at $p_0,\ldots,p_5$ are no worse than stated; there is an ordinary double point not in $\{p_0,\ldots,p_5\}$.

Let $T$ be a general line through $p_0.$ We have
$$
\begin{array}{l}
L_1:=8\widetilde{T}-4E_0-(2E_1+2E_1')-\sum_2^4(2E_i+3E_i')-2E_5,\\
L_2:=5\widetilde{T}-3E_0-\sum_1^3(E_i+2E_i')-(E_4+E_4')-E_5,\\
L_3:=4\widetilde{T}-2E_0-(E_1+2E_1')-\sum_2^3(E_i+E_i')-(E_4+2E_4'),
\end{array}
$$

$$
\begin{array}{l}
K_X+L_1\equiv \left(\widetilde{T_2}+\widetilde{T_3}+\widetilde{T_4}-3E_0-\sum_2^4 2E_i' \right) + \left(2\widetilde{T}-(E_1+E_1')-\sum_2^5 E_i\right),\\
K_X+L_2\equiv 2\widetilde{T}-2E_0-\sum_1^3 E_i',\\
K_X+L_3\equiv \widetilde{T}-E_0-E_1'-E_4'+E_5
\end{array}
$$
and
$$2K_X+\sum_1^3 L_i\equiv 11\widetilde{T}-7E_0-\sum_1^4(2E_i+4E_i')-E_5.$$

The divisor $$\widetilde{T_2}+\widetilde{T_3}+\widetilde{T_4}-3E_0-\sum_2^4 2E_i'$$ is a disjoint union of 6 $(-2)$-curves,
each meeting $K_X+L_1$ with intersection number $-1.$
Hence $K_X+L_1$ is effective only if $$2\widetilde{T}-(E_1+E_1')-\sum_2^5 E_i$$ is effective.
This is not the case, we can verify that the conic through the points $p_1,\ldots,p_5$ is not tangent to
the line $T_1.$ Therefore $h^0(X,\mathcal O_X(K_X+L_1))=0$ and then
$$p_g(V)=0+0+0+0=0.$$
Also
$$\chi(\mc O_V)=4+\frac{1}{2}(-2-2-2)=1$$
and
$$K_V^2=-9.$$
Since the minimal model $V'$ of $V$ is obtained contracting the $16$ $(-1)$-curves contained in
$\pi^*\left(\widetilde{T_1}+\cdots+\widetilde{T_4}\right),$ then $K_{V'}^2=7.$

The covering $\pi$ factors as $$V\longrightarrow Y\longrightarrow X,$$
where $Y\rightarrow X$ is the double cover with branch locus $D_2+D_3.$
Using the double cover formulas, one can verify that the smooth minimal model of $Y$ is a numerical
Campedelli surface ($p_g=q=0,$ $K^2=2$). The double cover $V\rightarrow Y$ is ramified over the
pullback of $D_1$ (which contains four $(-2)$-curves) and over the node corresponding to the
transverse intersection of $D_2$ and $D_3$.

\appendix
\section*{Appendix: Magma code}\label{appendix}

Here we use the computer algebra system Magma \cite{BCP} to show that the curves $C_6$ and $C_7$ referred in
Section \ref{pgq0} are exactly as stated there. This code can be tested on the online Magma calculator \cite{MC}.

\begin{small}

\begin{verbatim}
R<i>:=PolynomialRing(Rationals());
K<i>:=ext<Rationals()|i^2+1>;
P<x,y,z>:=ProjectiveSpace(K,2);

F6:=4*x^6-273*x^4*y^2-258*x^2*y^4-481*y^6+720*x^4*y*z+1740*x^2*y^3*z+
4020*y^5*z-520*x^4*z^2-3190*x^2*y^2*z^2-12670*y^4*z^2+1200*x^2*y*z^3+
17700*y^3*z^3+900*x^2*z^4-9225*y^2*z^4;

F7:=12*x^7+(8*i+420)*x^6*y+1611*x^5*y^2+(174*i+3060)*x^4*y^3+
4086*x^3*y^4+(924*i+3360)*x^2*y^5+987*x*y^6+(-242*i+720)*y^7-560*x^6*z-
4320*x^5*y*z+(-480*i-13580)*x^4*y^2*z-23940*x^3*y^3*z+
(-5160*i-24980)*x^2*y^4*z-10620*x*y^5*z+(1320*i-6960)*y^6*z+
2760*x^5*z^2+(240*i+16200)*x^4*y*z^2+44970*x^3*y^2*z^2+
(9780*i+63900)*x^2*y^3*z^2+39210*x*y^4*z^2+(-2460*i+25200)*y^5*z^2-
4400*x^4*z^3-28800*x^3*y*z^3+(-7200*i-62300)*x^2*y^2*z^3-
60300*x*y^3*z^3+(1800*i-40400)*y^4*z^3+2700*x^3*z^4+
(1800*i+16500)*x^2*y*z^4+33075*x*y^2*z^4+(-450*i+24000)*y^3*z^4;

C6:=Curve(P,F6); C7:=Curve(P,F7);
IsAbsolutelyIrreducible(C6);
IsAbsolutelyIrreducible(C7);

p:=[P![0,0,1],P![-2,1,1],P![2,1,1],P![-1,2,1],P![1,2,1],P![3,2*i,1]];

[ResolutionGraph(C6,p[i]):i in [1..5]];
[ResolutionGraph(C7,p[i]):i in [1..6]];
[ResolutionGraph(C6 join C7,p[i]):i in [1..6]];
SingularPoints(C6 join C7);
\end{verbatim}

\end{small}

\bibliography{ReferencesRito}

\

\

\noindent Carlos Rito\\
\\{\it Permanent address:}
\\ Universidade de Tr\'as-os-Montes e Alto Douro, UTAD
\\ Quinta de Prados
\\ 5000-801 Vila Real, Portugal
\\ www.utad.pt
\\ crito@utad.pt\\
\\{\it Current address:}
\\ Departamento de Matem\' atica
\\ Faculdade de Ci\^encias da Universidade do Porto
\\ Rua do Campo Alegre 687
\\ 4169-007 Porto, Portugal
\\ www.fc.up.pt
\\ crito@fc.up.pt

\end{document}